\documentclass[letterpaper, 10pt, conference]{ieeeconf}

\IEEEoverridecommandlockouts   

\usepackage{xspace,amssymb,epsfig,syntonly}
\usepackage{epsfig,amsmath,color}
\usepackage{xspace,syntonly,empheq}
\usepackage{bbm}
\usepackage{mathtools}
\usepackage{cite}
\usepackage{graphicx}
\usepackage{algorithm,algorithmic}
\usepackage{verbatim}
\usepackage{amssymb}
\usepackage{circuitikz}
\usepackage{enumerate}

\newcommand{\cN}{{\cal N}}

\newcommand{\cT}{{\cal T}}

\newcommand{\cE}{{\cal E}}

\newcommand{\cX}{{\cal X}}

\newcommand{\bx}{{\bf x}}

\newcommand{\sfT}{\textsf{T}}

\DeclarePairedDelimiterX{\norm}[1]{\lVert}{\rVert}{#1}
\DeclarePairedDelimiter{\abs}{\lvert}{\rvert}

\usepackage{epstopdf}

\newtheorem{lemma}{Lemma}
\newtheorem{theorem}{Theorem}

\newtheorem{assumption}{Assumption}

\title{\LARGE \bf Inexact Online Proximal-gradient Method \\ for Time-varying Convex Optimization}

\author{Amirhossein Ajalloeian, Andrea Simonetto, Emiliano Dall'Anese
\thanks{A. Ajalloeian and E. Dall'Anese are with the Department of Electrical, Computer, and Energy Engineering, University of Colorado Boulder, Boulder, CO, USA emails: amirhossein.ajalloeian@colorado.edu, emiliano.dallanese@colorado.edu. A. Simonetto is with  IBM Research Ireland, Dublin, Ireland; email: andrea.simonetto@ibm.com.}
\thanks{The work of A. Ajalloeian and E. Dall'Anese was supported by NREL APUP UGA-0-41026-109 and by the National Science Foundation Grant 1941896.}
}

\begin{document}

\maketitle

\thispagestyle{empty} 

\begin{abstract}
This paper considers an online proximal-gradient method to track the minimizers of a composite convex function that may continuously evolve over time. The online proximal-gradient method  is ``inexact,'' in the sense that: (i) it relies on an approximate first-order information of the smooth component of the cost; and, (ii)~the proximal operator (with respect to the non-smooth term) may be computed only up to a certain precision. Under suitable  assumptions, convergence of the error iterates is established for strongly convex cost functions. On the other hand, the dynamic regret is investigated when the cost is not strongly convex, under the additional assumption that the problem includes feasibility sets that are compact. Bounds are expressed in terms of the cumulative error and the path length of the optimal solutions. This suggests how to allocate resources to strike a balance between performance and  precision  in the gradient computation and in the proximal operator. 
\end{abstract}

\section{Introduction and Problem Formulation}
\label{sec:introduction}

The proximal-gradient method is a powerful framework for solving optimization problems with  objectives that consist of a differentiable convex function
and a nonsmooth convex function~\cite{Beck_FOM,parikh2014proximal,rockafellar1976monotone}. By taking advantage of this composite structure, proximal-gradient
methods are known to exhibit the same convergence rates of the gradient method for
 smooth problems~\cite{cevher2014convex}. Accordingly, proximal-gradient methods can be leveraged to efficiently solve a number of problems that arise in the broad areas of, e.g., statistical learning, network optimization, and design of optimal controllers for distributed systems~\cite{dixit2019online,dhingra2018proximal,fardad2011sparsity, schmidt2011convergence}.

This paper investigates the design of \emph{online} proximal-gradient methods for composite convex functions that continuously evolves over time. To outline the  setting concretely, discretize the temporal axis as $\cT := \{n \Delta, n \in \mathbb{N}\}$, with $\Delta$ a given interval. Let  $g_k: \cX \rightarrow  \mathbb{R}$ be a closed, convex and proper function with a Lipschitz-continuous gradient at each time $k \in \cT$, with $\cX \subseteq \mathbb{R}^n$ a given set; further, let $h_k:  \cX \rightarrow  \mathbb{R} \cup \{+ \infty\}$ be a lower semi-continuous proper convex function for all $k \in \cT$. Consider then the following time-varying optimization problem~\cite{popkov2005gradient,SimonettoGlobalsip2014,Simonetto17}:
\begin{equation}
\label{eq:problem_time_varying}
    \min_{x \in \mathbb{R}^n} f_k(x) := g_k(x) + h_k(x) \,, \hspace{.3cm} k \in \cT \, .
\end{equation}
Let $x_k^*$ be an optimal solution of~\eqref{eq:problem_time_varying} at time $k$  (which is unique if $f_k$ is strongly convex). In principle, a proximal method or an accelerated proximal method can be utilized to attain $x_k^*$; for example, it is known that when 
$f_t$ is convex $L$-smooth,
the number of iterations required to obtain an objective function within an error  $\xi$ is $\mathcal{O}(L \|x_{k,0} - x_k^*\|^2 / \xi )$ and $\mathcal{O}(\sqrt{L \|x_{k,0} - x_k^*\|^2 / \xi} )$ for a proximal method and its accelerated counterpart, respectively, with $x_{k,0}$ the starting point for the algorithm~\cite{cevher2014convex}\footnote{
\emph{Notation}: For a given vector $x \in \mathbb{R}^n$,  $\|\bx\| := \sqrt{\bx^\sfT \bx}$, with $^\sfT$ denoting transposition; for $x \in \mathbb{R}^n$ and $y \in \mathbb{R}^n$, $\langle x, y \rangle$ denotes the inner product.  For a differentiable function $f: \mathbb{R}^n \rightarrow \mathbb{R}$, $\nabla_{x} f(x)$ is the gradient vector of $f(x)$ with respect to $x \in \mathbb{R}^n$. If $f$ is non-differentiable, $\partial f(x)$ denotes the subdifferential of $f$ at $x$; in particular, a vector $v$ is a subgradient of $f$ at $x$ if $f(x^\prime) \geq f(x) + \langle v, x^\prime - x \rangle$ for all $x^\prime$ in the domain of $f$. On the other hand, $\partial_\epsilon  f(x)$ denotes the $\epsilon$-subdifferential of $f$ at $x$; a vector $v$ is an $\epsilon$-subgradient of $f$ at $x$ if $f(x^\prime) \geq f(x) + \langle v, x^\prime - x \rangle - \epsilon$ for all $x^\prime$ in the domain of $f$. A function $f: \mathbb{R}^n \rightarrow \mathbb{R}$ is $\mu$-strongly convex if $f(y) \geq f(x) + s_{x}^{\sfT}(y-x) + \frac{\mu}{2} \norm{y-x}^{2}$ for all $x, y$ and any $s_{x} \in \partial f(x)$. Finally, $\mathcal{O}$ refers to the big O  notation.
}. Results for strongly convex functions can be found in e.g.,~\cite{cevher2014convex,Beck_FOM}.

In contrast, this paper targets an online (or ``running'', or ``catching-up'' \cite{Moreau1977}) case where only one or a few steps of the proximal-gradient method  can be performed within an interval $\Delta$ (i.e., before the underlying optimization problem may change); further, the paper considers the case when  the implementation of the algorithmic steps is \emph{inexact}. Taking the case where only one step can be performed within an interval $\Delta$, an online  proximal-gradient algorithm amounts to the  execution of the following  steps~\cite{Simonetto17} at each time $k \in \cT$ :
\begin{subequations}
\label{eq:proximal_exact}
\begin{align}
    y_{k} & = x_{k-1} - \alpha \nabla_x g_{k}(x_{k-1}) \\
    x_{k} & =  \text{prox}_{\lambda h_{k}} \big\{y_{k}\big\}
\end{align}
\end{subequations}
where $\alpha > 0$ is the step size, and the proximal operator $\text{prox}_{\lambda h}: \mathbb{R}^n \rightarrow  \mathbb{R}^n$ is defined as~\cite{Beck_FOM}:
\begin{align}
\label{eq:proximal}
    \hspace{-.2cm} \text{prox}_{\lambda h} \big\{ y \big\} := \arg \min_{x} \Phi_{\lambda h}(x) := h(x) + \frac{1}{2 \lambda} \|x - y\|^2 \hspace{-.2cm}
\end{align}
with $\lambda > 0$ a given parameter. Notice that constraints of the form $x \in \cX_k$, with $\cX_k \subseteq \cX$ a convex  set, can be handled  via indicator functions~\cite{Beck_FOM}; i.e., by setting $h(x) = h^\prime(x) + \delta_{\cX_k}(x)$, where $h^\prime(x)$ is a lower-semicontinuous convex function, and $\delta_{\cX_k}(x) = 0$ if $x \in \cX_k$ and $\delta_{\cX_k}(x) = +\infty$ otherwise.

Inexactness of the steps~\eqref{eq:proximal_exact} may emerge because of the following two aspects: (i) only an approximate first-order information of $g_k$ may be available~\cite{Bernstein2019feedback,dixit2019online}; and, (ii) the proximal operator  may  be computed only up to a certain precision~\cite{schmidt2011convergence,salzo2012inexact,villa2013accelerated}. Before proceeding, examples of applications that motivate the proposed setting are briefly explained.

\vspace{.1cm}

\noindent \emph{Example 1: Feedback-based network optimization.} Online algorithms are, in this case, utilized to produce decisions to nodes of a networked system (i.e., a power system, a transportation network, or a communication network); temporal variability emerge from time-varying problem inputs (i.e., non-controllable power injections in a power system) or time-varying engineering objectives~\cite{Bernstein2019feedback,vaquero2018distributed}. Measurements of the network state are utilized to obtain an estimate of the gradient of $g_k$ at each time step. Inexactness of the proximal operator captures the case where the projection is performed onto an inner approximation of the actual feasibility region $\cX_k$ (i.e., when one has an approximate region for aggregations of energy resources)~\cite{nazir2018inner}; or, when the proximal operator is not easy to compute within an interval $\Delta$.

\vspace{.1cm}

\noindent \emph{Example 2: Online zeroth-order methods.} Zeroth-order methods involves an estimate of $\nabla_{x} g_k$ at $x$ based on functional evaluations $\{g_k(x+u_i), i = 1, \ldots, I\}$, with $u_i$ a given perturbation; see, e.g., Gaussian smoothing or  Kiefer-Wolfowitz approaches~\cite{Hajinezhad19,sahu2018distributed,chen2018bandit}. Inexactness of~\eqref{eq:proximal} is due to projections onto a restriction of $\cX_k$~\cite{chen2018bandit} or when  ~\eqref{eq:proximal} is not solved to convergence.

\vspace{.1cm}

\noindent \emph{Example 3: Learning under information streams.} For applications with continuous streams of data, the interval $\Delta$ may coincide with the inter-arrival time of data points; because of an underlying limited  computational budget (compared to $\Delta$), one may afford one step of the proximal-gradient method and a limited number of algorithmic steps to solve~\eqref{eq:proximal}. Examples include  
singular value decomposition (SVD) based proxies~\cite{dixit2019online} or structured sparsity~\cite{jenatton2010proximal}.

\vspace{.1cm}

For static optimization settings, convergence of inexact proximal-gradient methods  has been investigated in, e.g.,~\cite{salzo2012inexact,villa2013accelerated,schmidt2011convergence,machart2012optimal} (see also pertinent references therein); in particular,~\cite{schmidt2011convergence} showed that the inexact proximal-gradient method can achieve the same rate of convergence of the exact counterpart if the error sequence decreases at appropriate rates.  In an online setting,~\cite{dixit2019online}  investigated the convergence of the  proximal-gradient method with an approximate knowledge of $\nabla_x g_{k}$ (but with an exact implementation of the proximal operator); strongly convex cost functions were considered.  In this paper,  we analyze the convergence of the online inexact proximal-gradient method with errors in both the computation of $\nabla_x g_{k}$ and~\eqref{eq:proximal}. In particular: 

\noindent $\bullet$ The results of~\cite{dixit2019online} are generalized for the case of errors in the proximal operator, and for the case of costs that are not strongly-convex.     

\noindent $\bullet$ The analysis of, e.g.,~\cite{schmidt2011convergence,salzo2012inexact,villa2013accelerated} in the context of batch optimization is extended to a time-varying setting considered, with  the temporal variability of solution paths~\eqref{eq:problem_time_varying} playing  a key role in the convergence rates.

Under suitable  assumptions, convergence of the error iterates is established for strongly convex cost functions. On the other hand, convergence claims are established in terms of dynamic regret when the cost is not strongly convex, under the additional assumption that the feasibility sets are compact.  Bounds are expressed in terms of the cumulative error and the path length of the optimal solutions. The role of the errors is emphasized in the bounds, thus suggesting how to allocate computational resources to strike a balance between performance and  precision  in the gradient computation and in the proximal operator.

\section{Online Inexact Algorithm}
\label{sec:algorithm}

The models for the errors in the computation of the gradient of $g_k$ and of the proximal operator are described first, followed by the online inexact proximal-gradient method.  

\noindent \emph{Gradient error.} For a given point $x \in \cX$, the first-order information of $g_k$ is available in the form of $\tilde{\nabla} g_{k}(x) = \nabla_x g_{k}(x) + e_k$, with $e_k$ denoting the gradient error. The error sequence $\{e_k, k \in \cT\}$ is assumed to be bounded. 

\noindent \emph{Error in the proximal step.} A point $x$ is an approximation of $\text{prox}_{\lambda h_{k}} \big\{y\big\}$ with a precision $\epsilon \geq 0$ if~\cite{salzo2012inexact}:  
\begin{equation}
\label{eq:approx_def}
   0 \in \partial_{\frac{\epsilon^{2}}{2\lambda}} \Phi_{\lambda h_k} (x) \, .
\end{equation}
It is useful to notice that equation \eqref{eq:approx_def} implies that~\cite{salzo2012inexact,schmidt2011convergence}
\begin{align}
\label{approx_err}
    \frac{\epsilon^{2}}{2\lambda} \geq \Phi_{\lambda h_{k}}(x)-\Phi_{\lambda h_{k}}(\text{prox}_{\lambda h_k}(y))
\end{align}
where $\Phi_{\lambda  h_{k}}(\text{prox}_{\lambda h_k}(y))$ corresponds to the case where the proximal operator is computed exactly; furthermore, since $\Phi_{\lambda h_{k}}$ is a $\frac{1}{\lambda}$-strongly convex function, one has that: 
\begin{equation}
\label{approx_sc}
    \frac{1}{2\lambda}\norm{x-\text{prox}_{\lambda h_{k}}(y)}^{2}   \leq \Phi_{\lambda h_{k}}(x)-\Phi_{\lambda h_{k}}(\text{prox}_{\lambda h_{k}}(y)).
\end{equation}
Together, equations \eqref{approx_sc} and \eqref{approx_err}  imply that 
\begin{align}
    \norm{x-\text{prox}_{\lambda h_{k}}(y)} \leq \epsilon \, .
\end{align}
See also Appendix~\ref{sec:proof_inexact proximal}. More details on~\eqref{eq:approx_def}  will be provided shortly.

With these definitions in place, the online inexact algorithm is presented next, where the parameter $\lambda$ is set to $\lambda = \alpha$ as in e.g.,~\cite{Beck_FOM,schmidt2011convergence,cevher2014convex}.

\vspace{.2cm}

\hrule

\hrule

\vspace{.1cm}

\noindent \textbf{Online inexact proximal-gradient algorithm}

\vspace{.1cm}

\hrule 

\vspace{.1cm}

\noindent Initialize $x_0$, $\alpha$, and set $\lambda = \alpha$.

\noindent For each $k \in \cT$:
 
\textbf{[S1]} Obtain estimate of the gradient $\tilde{\nabla} g_{k}(x_{k-1})$ 

\textbf{[S2]} Perform the following updates:
\begin{subequations}
\label{eq:proximal_inexact}
\begin{align}
    y_k & = x_{k-1} - \alpha \tilde{\nabla} g_{k}(x_{k-1}) \\
    x_{k} & \approx_{\epsilon_k} \text{prox}_{\lambda h_{k}} \big\{y_{k}\big\}
\end{align}
\end{subequations}

\textbf{[S3]} Go to \textbf{[S1]}.

\vspace{.1cm}

\hrule

\hrule

\vspace{.2cm}

At each time step $k$, the algorithm is assumed to have the availability of: i) an estimate of $\nabla_x g_{k}(x_{k-1})$, and, ii) the function $h_k$. This is the case when, e.g., $g_k$ depends on data or its gradient is not available, while $h_k$ represents  regularization terms based on a prior on the optimal solution, or set indicator functions for  constraints. It is thus reasonable to assume that one has access to $h_k$ -- since it is in general engineered -- while access on $\nabla_x g_k$ depends on data; take for example, the case $\|x-b_k\|^2_2 + \|x\|_1$, where $b_k$ is the data stream. Conditions on the step size $\alpha$ will be given shortly in Section~\ref{sec:performance}.

The characteristics  of error sequence $\{e_k, k \in \cT\}$ depends on the particular application. For example, in measurements-based online network optimization algorithms, $e_k$ captures measurement noise (see Example~1)~\cite{Bernstein2019feedback,vaquero2018distributed}; therefore, a bound on $\|e_k\|$~\cite{Bernstein2019feedback} (or on the expected value of $\|e_k\|$~\cite{dixit2019online}) is utilized to assess the tracking performance of the algorithm, but it may not be under the control of the designer of the algorithm. On the other hand, the error $e_k$ may be controllable by the designer of the algorithm in, e.g., zeroth-order methods (Example 2) and applications such as subspace  tracking and online sparse regression (Example 3); see e.g.,~\cite{Hajinezhad19,sahu2018distributed,chen2018bandit,dixit2019online,jenatton2010proximal} and pertinent references therein. 

Regarding the error sequence $\{\epsilon_k\}$, there are two common themes in the examples considered in Section~\ref{sec:introduction}: i) if a set indicator function $\delta_{\cX_k}(x)$ is considered, points may be projected in the interior of $\cX_k$; and, ii) for a given lower semi-continuous convex function,~\eqref{eq:proximal} may not be solved to convergence. In both cases, the error sequence $\{\epsilon_k, k \in \cT\}$ can be controlled, based on given computational budgets or other design specifications.  Examples are provided next. 

\vspace{.1cm}

\noindent \emph{Example: inexact projection}~\cite{salzo2012inexact}. Suppose that $h_k = \delta_{\cX_k}$, for a given closed and convex set $\cX_k$; let $d(y, \cX)$ denote the distance of the point $y$ from the convex set $\cX$. Then, the definition~\eqref{eq:approx_def} implies that $x \approx_{\epsilon} \text{prox}_{\lambda h}(y)$ if and only if 
\begin{align}
    x \in \cX_k ~\textrm{~and~} \|x - y\|^2 \leq d(y, \cX)^2 + \epsilon^2 \, .
\end{align}
That is, if $y \notin \cX_k$, then $x$ may not lie in the boundary of $\cX_k$; rather, $x$ may lie in the interior of $\cX_k$. It is worth noticing that the point $x$ is always feasible as explained in~\cite{salzo2012inexact}.

\vspace{.1cm}

\noindent \emph{Example: structured sparsity}. Take, for example, the case where $h_k(x) = \sum_i \|[x_{i_1}, \ldots, x_{i_I}]^\sfT\|_2$ with $[x_{i_1}, \ldots, x_{i_I}]^\sfT$ a given sub-block of the vector $x$. In this case, a block coordinate method can be utilized to solve~\eqref{eq:proximal_exact}~\cite{schmidt2011convergence,jenatton2010proximal}; the block coordinate method can be run up a given  error $\epsilon_k$.  

\vspace{.1cm}

\noindent \emph{Example: SVD-based proxies}. Proximal operators that involve an SVD computation (e.g., nuclear norm minimization) may be computed inexactly, especially for large matrices.  

\vspace{.1cm}

\noindent \emph{Example: distributed computation analyzed as an inexact centralized method}.  An additional motivation for considering inexact steps in the algorithm emerges from~\cite{Fast_distributed_PGM}; in particular,~\cite{Fast_distributed_PGM} considers a distributed proximal-gradient method for minimizing the sum of the cost functions of individual agents in a network. At each iteration of the algorithm, each agent updates its estimate along the negative gradient of the differentiable part of its cost function; then, a consensus step is performed, followed  by a local proximal step with respect to the nondifferentiable part of the local cost functions. Indeed, as discussed in~\cite{Fast_distributed_PGM}, this setup can be seen as an inexact centralized proximal gradient algorithm where the error sequences in the gradient and the proximal operator depend on the accuracy of the consensus step.

\vspace{.1cm}

Finally, it is also worth mentioning that the two errors could be analyzed in a unified way if one interprets $e_k$ as a perturbation in the computation of the (exact) operator $\text{prox}_{\lambda h_{k}}(x_{k-1} - \alpha \nabla_x g_{k}(x_{k-1}) )$ (see, e.g., Definition 3 in~\cite{salzo2012inexact}); however, similarly to~\cite{schmidt2011convergence}, the current models of the errors allows one to better appreciate the role of the ``exactness'' of the first-order information and the proximal operator in the performance of the algorithm as shown in the next section. 

\section{Performance Analysis}
\label{sec:performance}

This section will analyze the performance of the online inexact algorithm~\eqref{eq:proximal_inexact}; two metrics will be considered: 

\noindent \emph{i)} convergence of the  sequence $\{\|x_k - x_k^*\|, k \in \cT\}$; and,

\noindent \emph{ii)} the  dynamic regret, defined as (see, e.g.,~\cite{hall2015online,yi2016tracking,Jadbabaie2015} and references therein):
\begin{align}
\label{eq:regret}
\textrm{Reg}_k := \sum_{i = 1}^k f_{i}(x_{i}) - f_{i}(x_i^*)  \, .  
\end{align}
The error sequence $\{\|x_k - x_k^*\|, k \in \cT\}$ will be analyzed when the cost function is strongly convex; in particular, bounds on the cumulative error $\sum_{i = 1}^k \|x_i - x_i^*\|$ and $Q$-linear convergence results will be offered. When the cost function is not strongly convex and pertinent relaxed conditions for linear convergence (see, e.g.,~\cite{necoara2019linear} for a quadratic functional
growth condition) are not satisfied, $Q$-linear convergence may not be available. In that case, the dynamic regret can be used as a performance metric. 

The following standard assumptions are presumed throughout this section.  

\vspace{.1cm}

\begin{assumption}
\label{as:function_g} 
The function $g_k: \cX \rightarrow  \mathbb{R}$  is closed, convex and proper. Assume that $g_k$ 
has a $L_k$-Lipschitz continuous gradient at each time $k \in \cT$, and there exists $L$ such that $L_k < L$  for all $k \in \cT$.
\end{assumption}

\vspace{.1cm}

\begin{assumption}
\label{as:function_h} 
The function $h_k: \cX \rightarrow  \mathbb{R} \cup \{+ \infty\}$ is a lower semi-continuous proper convex function for all $k \in \cT$.
\end{assumption}

\vspace{.1cm}

\begin{assumption} 
\label{as:minimum_attained}
For  all $k \in \cT$, $f_{k}=g_{k}+h_{k}$ attains its minimum at a certain $x_{k}^{*} \in \cX_k$.
\end{assumption}

\vspace{.1cm}

To characterize  bounds on the error sequence and the dynamic regret, it is necessary to introduce a ``measure'' of the temporal variability of~\eqref{eq:problem_time_varying} as well as  of the ``exactness'' of the first-order information of $\{g_k, k \in \cT\}$ and the computation of the proximal operator. For the former, define   
\begin{align}
\sigma_k := \norm{x_{k}^{*}-x_{k-1}^{*}}
\end{align}
along with the following quantities~\cite{dixit2019online,Bernstein2019feedback,yi2016tracking}: 
\begin{align}
\label{eq:path}
\Sigma_k := \sum_{i = 1}^k \sigma_i , \, \hspace{.2cm} \bar{\Sigma}_k := \sum_{i = 1}^k \sigma_i^2
\end{align}
with $\Sigma_k$ typically referred to as the ``path length'' or ``cumulative drifting.'' When $f_k$ is strongly convex, $\sigma_k$ is uniquely defined; on the other hand, $\sigma_k$ is associated with a solution path when $f_k$ is not strongly convex. Consider further the following definitions for the cumulative errors~\cite{schmidt2011convergence}: 
\begin{align}
\label{eq:cumulative_error}
E_k := \sum_{i = 1}^k \norm{e_{i}}  , \hspace{.2cm}
P_k := \sum_{i = 1}^k \epsilon_i , \, \hspace{.2cm} \bar{P}_k := \sum_{i = 1}^k \epsilon_i^2  \, . \, 
\end{align}

With these definitions in place, the convergence results are established first for the case where the function $f_k$ in~\eqref{eq:problem_time_varying}  is $\mu_k$-strongly convex for all $k \in \cT$. 

The following lemma will be utilized to derive convergence results when the function $f_k$ in~\eqref{eq:problem_time_varying}  is $\mu_k$-strongly convex.

\vspace{0.2cm}

\begin{lemma}
\label{lemma:strongly_convex}
Let  Assumptions~\ref{as:function_g}--\ref{as:minimum_attained} hold, and assume that $g_k$ is $\mu_k$-strongly convex for all $k \in \cT$. Then, the following holds for the algorithm~\eqref{eq:proximal_inexact}: 
\begin{equation}
    \norm{x_{k}-x_{k}^{*}} \leq \rho_{k} \norm{x_{k-1}-x_{k-1}^{*}} + \rho_{k} \sigma_{k} + \alpha \norm{e_{k}} + \epsilon_{k}
    \label{before_recur1}
\end{equation}
where in \eqref{before_recur1}, $\rho_k := \max\{|1-\alpha \mu_k|, |1-\alpha L_k|\}$. Applying \eqref{before_recur1} recursively, it holds that:
\begin{equation}
\label{after_recur1-1}
    \norm{x_{k}-x_{k}^{*}} \leq  \beta_k \norm{x_{0}-x_{0}^{*}} +  \sum_{i=1}^{k} \eta_{k,i} \sigma_{i} + \sum_{i=1}^{k} \nu_{k,i}(\alpha \norm{e_i} + \epsilon_i) 
\end{equation}
where in \eqref{after_recur1-1}, $\beta_k := \Pi_{i = 1}^{k} \rho_{i}$, $\eta_{k,i}:= \prod_{\ell = i}^k \rho_{\ell}$ and
\begin{align}
\nu_{k,i} := 
\begin{cases}
\prod_{\ell = i + 1}^k \rho_{\ell}, & \textrm{if~} i = 1, \ldots, k - 1 \\
1, &  \textrm{if~} i = k , 
\end{cases}
\end{align} 
for all $i \in \cT$.
\end{lemma}
\emph{Proof}. See Appendix~\ref{sec:proof_lemma_strongly_convex}.

\vspace{.2cm}

Based on Lemma~\ref{thm:strongly_convex}, the following theorems characterize the behavior of the error sequence $\{ \norm{x_{i}-x_{i}^{*}}, k \in \cT\}$. 

\vspace{.2cm}

\begin{theorem}
\label{thm:strongly_convex}
Suppose that Assumptions~\ref{as:function_g}--\ref{as:minimum_attained} hold. Assume that $g_k$ is $\mu_k$-strongly convex for all $k \in \cT$, and that $\mu < \mu_k$ for all $k$, for a given $\mu > 0$.  If $\alpha < 2 / L$, then 
\begin{align}
\label{eq:strongly_convex}
& \hspace{-.3cm} \sum_{i = 1}^k  \norm{x_{i}-x_{i}^{*}} \leq \frac{1}{1 - \rho} \left(\rho \norm{x_{0}-x_{0}^{*}} + \rho \Sigma_k + P_k + \alpha E_k \right)  
\end{align}
where $\rho := \sup_k \{\rho_k\} < 1$. 
\end{theorem}
\vspace{.1cm}
\emph{Proof}. See Appendix~\ref{sec:proof_stronglyconvex}.

\vspace{.2cm}

\begin{theorem}
\label{thm:strongly_convex_maximum}
Suppose that there exists finite constants $\sigma$, $\gamma_{e}$ and $\gamma_{\epsilon}$ such that $\sigma_k \leq \sigma$, $\norm{e_k} \leq \gamma_{e}$, and $\epsilon_k \leq \gamma_{\epsilon}$ for all $k \in \cT$. Then, under the same assumptions of Theorem~\ref{thm:strongly_convex}, it hols that
\begin{equation}
\label{eq:strongly_convex2}
    \limsup_{k \to \infty} \norm{x_{k}-x_{k}^{*}} \leq \frac{\alpha \gamma_{e} + \gamma_{\epsilon} + \rho \sigma}{1-\rho} . 
\end{equation}
\end{theorem}
\vspace{.1cm}
\emph{Proof}. See Appendix~\ref{sec:proof_stronglyconvex_maximum}.

\vspace{.2cm}

From Theorem~\ref{thm:strongly_convex}, it can be seen that if $\Sigma_k$, $P_k$, and $E_k$ grow as $\mathcal{O}(k)$, then the averaged tracking error behaves as $(1/k)\sum_{i = 1}^k  \norm{x_{i}-x_{i}^{*}} = \mathcal{O}(1)$. The same limiting behavior can be obtained even if the error sequences $\{e_k\}$ and $\{\epsilon_k\}$ decrease over time, if $\Sigma_k = \mathcal{O}(k)$.  As expected, the error in the gradient computation is down-weighted by the step size; on the other hand, the error in the proximal operator directly affects the tracking performance. The result of Theorem~\ref{thm:strongly_convex_maximum} may suggest how to allocate computational resources to minimize the maximum tracking error, for a given $\alpha$ and $\rho$; since $\sigma$ is multiplied by $\rho$, one may want to increase the interval $\Delta$ (thus increasing $\sigma$) and allocate more resources in the proximity operator (thus decreasing $\gamma_\epsilon$).

The next result pertains to the dynamic regret, and it extends the existing results of~\cite{dixit2019online} to the case of inexact proximal operators. 

\vspace{.2cm}

\begin{theorem}
\label{thm:regret_strongly_convex}
Suppose that Assumptions~\ref{as:function_g}--\ref{as:minimum_attained} hold  and assume that there exists $D_k < + \infty$ such that $\|\partial f_k\| \leq D_k$ over $\cX$. Then, if $f_{k}$ is strongly convex for all $k\in \cT$, the dynamic regret has the following limiting behavior
\begin{align}
\label{eq:regret_strongly_convex}
\textrm{Reg}_k = \mathcal{O}\left(1 + \Sigma_k + P_k + E_k \right)   \, .
\end{align}
\end{theorem}
\vspace{.1cm}
\emph{Proof}. See Appendix~\ref{sec:proof_regret_strongly_convex}.

\vspace{.2cm}

Finally, the next result pertains to the dynamic regret in case of functions $f_k$ that are convex but not strongly-convex. We impose the additional assumption that the cost function includes a time-varying set indicator function for a compact set. The derivation of bounds on the dynamic regret for the case of sets that are not compact is left as a future research. 

\vspace{.2cm}

\begin{theorem}
\label{thm:regret_convex}
Suppose that Assumptions~\ref{as:function_g}--\ref{as:minimum_attained} hold, and suppose that $f_{k}$ is convex and $h(x) = h^\prime(x) + \delta_{\cX_k}(x)$, where $h^\prime(x)$ is a lower-semicontinuous convex function and $\cX_k \subseteq \cX$ is compact for all $k \in \cT$. Suppose that $\alpha \leq 1/\sup\{L_k\}$. Then, the dynamic regret has the following limiting behavior
\begin{align}
\label{eq:regret_strongly_convex}
\textrm{Reg}_k = \mathcal{O}\left(1 + k + \Sigma_k + \bar{\Sigma}_k + P_k + \bar{P}_k + E_k \right)   \, .
\end{align}
\end{theorem}
\vspace{.1cm}
\emph{Proof}. See Appendix~\ref{sec:proof_regret_convex}.

\vspace{.2cm}

If $\Sigma_k$, $\bar{\Sigma}_k$, $P_k$, $\bar{P}_k $, and $E_k$ grow as $\mathcal{O}(k)$, then $(1/k)\textrm{Reg}_k = \mathcal{O}(1)$; that is, the  dynamic regret settles to a constant value. A no-regret result can be obtained if $\Sigma_k$, $\bar{\Sigma}_k$, $P_k$, $\bar{P}_k $, and $E_k$ all grow sublinearly  in $k$; that is, if they grow as $o(k)$. However, this may not be achievable when the optimization problem continuously evolves over time.

\section{Illustrative Numerical Results}
\label{sec:results}

\begin{figure}[t!]
  \centering
  \includegraphics[width=.4\columnwidth]{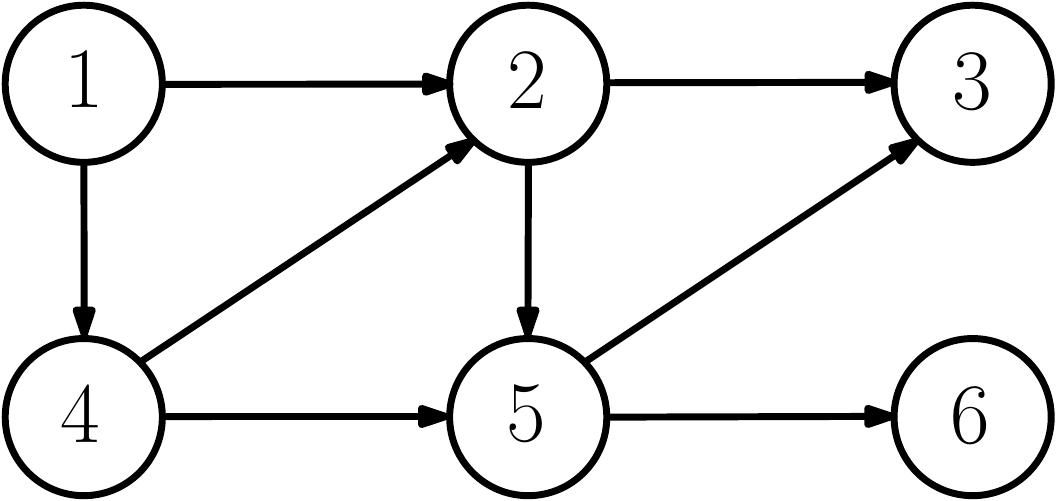}
\caption{Network utilized in the numerical results. Two data traffic flows are generated at nodes $1$ and $4$, with destination $3$ and $6$, respectively.}
\label{fig:F_comm_net}
\vspace{-.2cm}
\end{figure}

As an example of application of the proposed methods, we consider a network flow problem based on the network in Fig.~\ref{fig:F_comm_net}. the network graph $(\cN, \cE)$ the network has $|\cN| = 6$ nodes and $|\cE| = 8$ (directed) links, and the routing matrix $T$ is based on the directed edges. Let $z_k(i,s)$ denote the  rate generated at node $i$ for traffic $s$ at time $k \in \cT$, and $x_k(ij,s)$ the flow between noted $i$ and $j$ for traffic $s$. For brevity, let $z_k(s)$ and $x_k(s)$ stacks the node traffic $\{z_k(i,s), i \in \cN\}$ and link rates $\{x_k(ij,s), (i,j \in \cE)\}$ for the $s$-th flow. One has that $z_k(s)$ and $x_k(s)$ are related by the flow conservation constraint $z_k(s) = T x_k(s)$; if a node $i$ does not generate or receive traffic, then $z(i,s) = 0$ or $0 = T_i x(s)$ with $T_i$ the $i$-th row of $T$.  
Consider then the following time-varying problem (where we recall that $k$ is the time index):
\begin{align} 
\label{eqn:sampledProblemComm}
&\min_{\substack{\{z, x \} } }\sum_{i,s} - \kappa_k(i,s) \log(1+z(i,s)) + \frac{\nu}{2} (\|z\|_2^2 + \|x\|_2^2 ) \nonumber \\
& \hspace{.5cm} + \sum_{i,s} \delta_{\{0 \leq \sum_s  x(ij,s)  + w_k(ij) \leq  c_k(ij)\}} + \delta_{\{0 \leq z \leq z^{\textrm{max}}\}} \nonumber \\
& \hspace{.5cm} + \sum_s \delta_{\{z(s) = T x(s)\}}  
\end{align} 
where $z$ and $x$ stack the traffic rates and link rates for brevity; 
$z^{\textrm{max}}$ is a maximum traffic rate; and  $\nu > 0$ is a regularization function that makes the cost strongly convex. The per-link capacity constraints is $0 \leq \sum_{(i,j) \in \cE} x(ij,s) + w_k(ij) \leq c_k(ij)$, where $c_k(i,j)$ is the time-varying link capacity and  $w_k(ij)$ is a time-varying link traffic that is non-controllable. Notice that problem~\eqref{eqn:sampledProblemComm} can be equivalently rewritten in terms of only the vector variable $x$ (and the term $\|z\|_2^2$ in the cost can be dropped); simulations will be based on the reformulated problem.

For the numerical results,  assume that two traffic flows are generated by nodes $1$ and $4$, and they are received at nodes $3$ and $6$, respectively. In terms of dynamics of the optimal solutions, at each time step the channel gain of links are generated by using a complex Gaussian random variable with mean $1 + \jmath 1$ and a given variance $10^{-2}$ for both real and imaginary parts; the transmit power for each node is a  Gaussian random variable with mean $1$ and a variance $10^{-3}$; the exogenous traffic follows a random walk, where the increment has zero mean and a variance $10^{-2}$; and, the cost is perturbed by modifying $\kappa_k(i,s)$. Different values for $\sigma_t$ and $\sigma$ are obtained by varying the variance of these random variables.

\begin{figure}[t!]
  \centering
  \includegraphics[width=1.0\columnwidth]{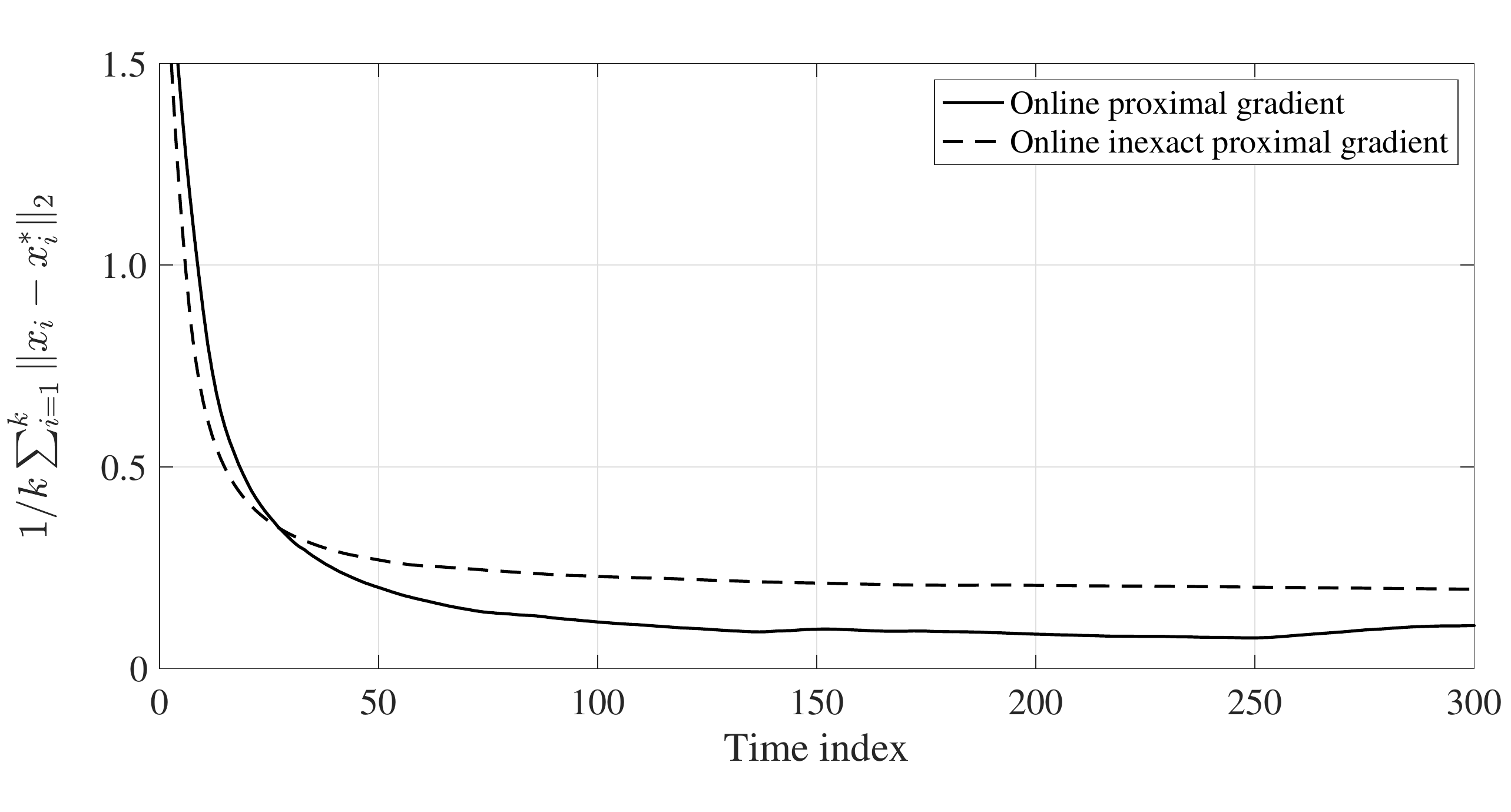}
\vspace{-.4cm}
\caption{Evolution of $(1/k) \sum_{i = 1}^k  \norm{x_{i}-x_{i}^{*}}$ for the case of exact online proximal gradient method and inexact proximal gradient method. }
\label{fig:F_acc2020_err_x}
\vspace{-.2cm}
\end{figure}

\begin{figure}[t!]
  \centering
  \includegraphics[width=1.0\columnwidth]{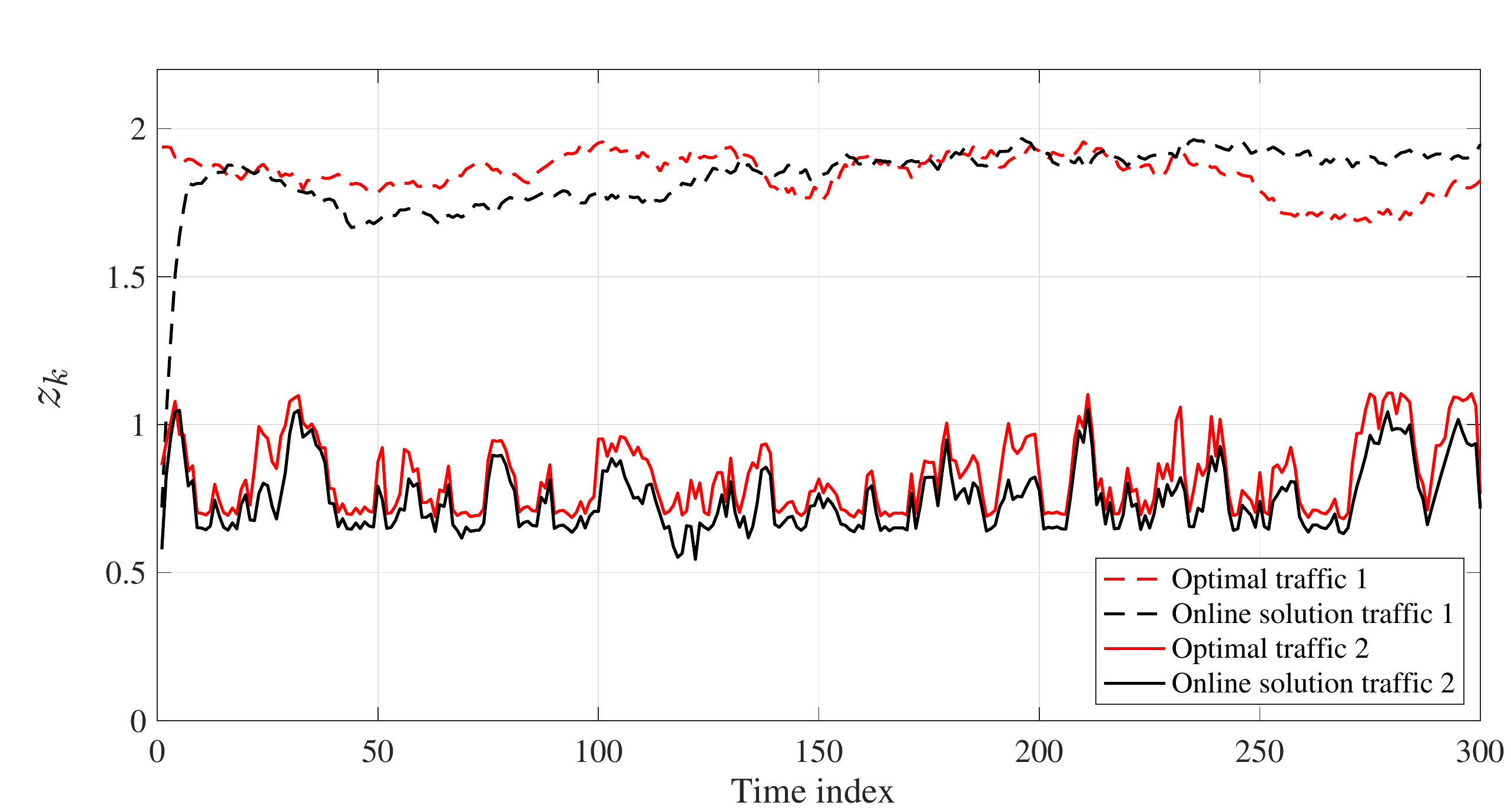}
 \vspace{-.2cm}
\caption{Example of traffic achieved for the case of batch solution (red) and online solution (black).}
\label{fig:F_acc2020_traffic_z}
\vspace{-.2cm}
\end{figure}

The algorithm~\eqref{eq:proximal_inexact} was implemented, with the following settings: 

\noindent $\bullet$ Gradient errors: the cost function $\kappa_k(i,s) \log(1+z(i,s))$ was assumed unknown; therefore, at each step of the algorithm, the gradient is estimated using a multi-point  bandit feedback~\cite{Hajinezhad19,chen2018bandit}. Briefly, to estimate the gradient of a function $g_k$ around a point $x_{k-1}$, consider drawing $M-1$ points $\{u_i\}_{i = 1}^{M-1}$ from the unit sphere; then, an estimate can be found as $ \tilde{\nabla} g_{k}(x_{k-1}) = (n/s(M-1)) \sum_{i = 1}^{M-1} (g_{k}(x_{k-1} + s u_i) - g_{k}(x_{k-1}))u_i$, where $s > 0$ is given parameter~\cite{chen2018bandit}. Notice that this requires    $M$ functional evaluations at each step of the algorithm. 

\noindent $\bullet$ Error in the proximal operator:  since the estimate of the gradient requires functional evaluations around the current point, we consider a restriction of the feasible set; in particular, consider the constraint $\rho \leq \sum_s  x(ij,s)  + w_k(ij) \leq  c_k(ij) - \rho$, where $\rho > 0$ is a pre-selected constant (that is related to $s$~\cite{chen2018bandit}). This will allow for functional evaluations $g_{k}(x_{k-1} + s u_i)$ at points $\{x_{k-1} + s u_i\}_{i = 1}^{M-1}$ that are feasible. In the numerical tests, $\sigma$ amounts to $0.7$, whereas the maximum errors due to the gradient estimate and the inexact projection add up to $0.5$.   

Figure~\ref{fig:F_acc2020_err_x} shows the evolution of the cumulative tracking error $(1/k) \sum_{i = 1}^k  \norm{x_{i}-x_{i}^{*}}$ at each time step $k$ (with $x_{k}^{*}$ unique, since the cost function is by design strongly convex). Based on Theorem~\ref{thm:strongly_convex}, in the current setting the limiting behavior of $(1/k) \sum_{i = 1}^k  \norm{x_{i}-x_{i}^{*}}$ is $\mathcal{O}(1)$. Indeed, a plateau can be seen, with an asymptotic error that is larger for the inexact proximal gradient method. Figure~\ref{fig:F_acc2020_traffic_z} illustrates the traffic rates achieved with a batch algorithm and with the inexact proximal gradient method.  It can be seen that the optimal traffic rates are tracked. Slightly lower traffic rates are obtained in the online case because of the projection onto a restriction of the feasible set. 

\section*{Acknowledgements}

The authors would like to thank the anonymous reviewers for the feedback provided on the paper.  

\appendix

\subsection{Results for the inexact proximal operator}
\label{sec:proof_inexact proximal}

Technical details for the error in the proximal operator are derived. These  technical details are  then utilized in Section~\ref{sec:performance} to derive pertinent convergence results. We start from the following lemma. 

\vspace{.2cm}

\begin{lemma}
\label{lem:epsilon-subdifferential}
Take the  function $f(x) = f_{1}(x)+f_{2}(x)$. Then $\partial_{\epsilon}f(x) \subset \partial_{\epsilon}f_{1}(x) + \partial_{\epsilon}f_{2}(x)$.

\textit{Proof}. The proof can be derived as a special case of Theorem 3.1.1 in \cite{hiriart1996convex}.
\end{lemma}

\vspace{.2cm}

Based on Lemma \ref{lem:epsilon-subdifferential}, we obtain:
\begin{multline}
   0 \in \partial_{\frac{\epsilon_k^{2}}{2\lambda}} \Phi_{\lambda h_k} (x_k) \subset \partial_{\frac{\epsilon_{k}^{2}}{2\lambda}} \big\{ \frac{1}{2\lambda}\norm{x_{k}-y_{k}}^{2} \big\} + \partial_{\frac{\epsilon_{k}^{2}}{2\lambda}}h_{k}(x_{k})  \\ \implies
 \hskip-1mm -\zeta \in \partial_{\frac{\epsilon_{k}^{2}}{2\lambda}}h_{k}(x_{k}), \quad  \zeta \in \partial_{\frac{\epsilon_{k}^{2}}{2\lambda}} \big\{ \frac{1}{2\lambda}\norm{x_{k}-y_{k}}^{2} \big\}.
\end{multline}
Since $\zeta$ is an $\frac{\epsilon_{k}^{2}}{2\lambda}$ subgradient:
\begin{equation}
    \frac{1}{2\lambda} \norm{x_{k}-y_{k}}^{2} - \frac{1}{2\lambda} \norm{t-y_{k}}^{2} + \zeta^{\sfT}(t-x_{k}) \leq \frac{\epsilon^{2}}{2\lambda} ~~~ \forall~t.
    \label{e_subgrad}
\end{equation}
Since \eqref{e_subgrad} is true for every  $t$, it is also specifically true for $t = y_{k} + \lambda \zeta$. Then we can write \eqref{e_subgrad} as:
\begin{equation}
    \frac{1}{2\lambda} \big( \norm{x_{k}-y_{k}}^{2} + \lambda^{2}\norm{\zeta}^{2} - 2\lambda\zeta^{\sfT}(x_{k}-y_{k})\big) \leq \frac{\epsilon_{k}^{2}}{2\lambda}
    \label{e_subgrad2}
\end{equation}
and therefore 
%
   $\norm{x_{k}-y_{k}-\lambda\zeta}^{2}/{2\lambda} \leq {\epsilon_{k}^{2}}/{2\lambda}$.
Now defining $r_{k}:= x_{k}-y_{k}-\lambda\zeta$ and putting \eqref{e_subgrad} and \eqref{e_subgrad2} together, we get
\begin{equation}
\label{approx_inclusion}
   \frac{y_{k}-x_{k}+r_{k}}{\lambda} \in  \partial_{\frac{\epsilon_{k}^{2}}{2\lambda}}h_{k}(x_{k}) \hspace{0.2cm} \textrm{with} \hspace{0.2cm} \norm{r_{k}} \leq \epsilon_{k}.
\end{equation}
Since $\lambda$ is chosen equal to $\alpha$ in \eqref{eq:proximal_inexact}, then applying the $\epsilon$-subdifferential definition for $h_{k}(x_{k})$, \eqref{approx_inclusion} can be written as the following inequality:
\begin{equation}
\label{for_NSC_proof}
\begin{array}{l}
    h_{k}(x_{k}) \leq h_{k}(x_{k}^{*}) - \frac{1}{\alpha} \langle x_{k-1}-x_{k}+r_{k},x_{k}^{*}-x_{k} \rangle 
    \\
    \hspace{2.2cm} + \langle \nabla g_{k}(x_{k-1})+e_{k},x_{k}^{*}-x_{k} \rangle + \frac{\epsilon_{k}^{2}}{2\alpha}
\end{array}
\end{equation}
where $\norm{r_{k}} \leq \epsilon_{k}$. If there is no inexactness, ${x_{k} = y_{k} + \alpha \zeta}$ is the exact proximal at $y_{k}$, i.e., $\text{prox}_{\alpha h_{k}}(y_{k})$. Therefore, the inexact proximal at $y_{k}$, i.e., $x_{k} \approx_{\epsilon} \text{prox}_{\alpha h_{k}}(y_{k})$ can be written as:
\begin{equation}
\label{for_SC_proof}
\begin{array}{l}
    x_{k} = \text{prox}_{\alpha h_{k}} (y_k) + r_{k} 
\end{array}
\end{equation}
with $\norm{r_{k}} \leq \epsilon_{k}$.

\subsection{Proof of Lemma~\ref{lemma:strongly_convex}}
\label{sec:proof_lemma_strongly_convex}
Based on \eqref{for_SC_proof}, we can write
\begin{subequations}
\begin{align}
    \nonumber \norm{x_{k}-x_{k}^{*}} & = \big\lVert\text{prox}_{\alpha h_{k}}[x_{k-1} - \alpha \tilde{\nabla}g_{k}(x_{k-1})] \\
    &\hspace{.5cm} + r_{k} - \text{prox}_{\alpha h_{k}}[x_{k}^{*} - \alpha \nabla g_{k}(x_{k}^{*})]\big\rVert
    \\
    &  \nonumber \leq \big \lVert \text{prox}_{\alpha h_{k}}[x_{k-1} - \alpha \tilde{\nabla}g_{k}(x_{k-1})] 
    \\
    &  \hspace{.5cm}- \text{prox}_{\alpha h_{k}}[x_{k}^{*} - \alpha \nabla g_{k}(x_{k}^{*})] \big \rVert + \epsilon_k \, .
\end{align}
\end{subequations}
Now, the prox operator is non-expansive, therefore:
\begin{multline}
    \norm{x_{k}-x_{k}^{*}} \leq \\ \big \lVert [x_{k-1} - \alpha \tilde{\nabla}g_{k}(x_{k-1})] - [x_{k}^{*} - \alpha \nabla g_{k}(x_{k}^{*})]  \big \rVert + \epsilon_k,
    \label{eq:dummy}
\end{multline}
which leads to 
\begin{multline}
    \norm{x_{k}-x_{k}^{*}} \leq  \alpha \|e_k\| + \epsilon_k + \\ \big \lVert [x_{k-1} - \alpha {\nabla}g_{k}(x_{k-1})] - [x_{k}^{*} - \alpha \nabla g_{k}(x_{k}^{*})]  \big \rVert.
    \label{eq:dummy1}
\end{multline}
Consider the function $\xi(x) := x - \alpha \nabla g_k(x)$. The norm of the gradient $\nabla \xi(x)$ is bounded as
\begin{equation}
\|\nabla \xi(x)\| \leq \max\{|1-\alpha \mu_k|, |1-\alpha L_k|\} =: \rho_k,    
\end{equation}
and therefore $\xi(x)$ is Lipschitz (and a contraction for $\alpha < 2/L_k$)\cite{ryu2016primer}. Hence, we can bound~\eqref{eq:dummy1} as
\begin{equation}
    \norm{x_{k}-x_{k}^{*}} \leq  \alpha \|e_k\| + \epsilon_k + \rho_k \norm{x_{k-1} - x_{k}^{*}}.
    \label{eq:dummy2}
\end{equation}

Adding and subtracting $x_{k-1}^{*}$ to the first term in the right hand side of \eqref{eq:dummy2} we can rewrite \eqref{eq:dummy2} as
\begin{subequations}
\begin{align}
    & \norm{x_{k}-x_{k}^{*}} \leq \rho_{k} \norm{x_{k-1}-x_{k-1}^{*}+x_{k-1}^{*}-x_{k}^{*}} + \alpha \norm{e_{k}} + \epsilon_{k}
    \\
    & \hspace{0.5cm} \leq \rho_{k} \norm{x_{k-1}-x_{k-1}^{*}} + \rho_{k} \norm{x_{k-1}^{*}-x_{k}^{*}}+ \alpha \norm{e_{k}} + \epsilon_{k}
    \\
    & \hspace{0.5cm} \leq \rho_{k} \norm{x_{k-1}-x_{k-1}^{*}} + \rho_{k} \sigma_{k} + \alpha \norm{e_{k}} + \epsilon_{k}.  \label{recursion1}
\end{align}
\end{subequations}

Applying \eqref{recursion1} recursively we get \eqref{after_recur1-1}.

\subsection{Proof of Theorem~\ref{thm:strongly_convex}}
\label{sec:proof_stronglyconvex}
Since $\alpha < 2 / L$ for all $k \in \cT$, one has that  $\rho_{k}<1$ for all $k \in \cT$. Also, one can define $\rho := \sup_k \{\rho_k\} < 1$. It is always possible to define such a $\rho$, because $L_{k} < L$ and $\mu_{k} > \mu$ for all $k$, and therefore both $\alpha L_{k}$ and $\alpha \mu_{k}$ are upper bounded by $\alpha L$ and lower bounded by $\alpha \mu$. As $\alpha L$ and $\alpha \mu$ are strictly positive, therefore we can always define $\rho := \sup_k \{\rho_k\} < 1$. Summing both sides of \eqref{before_recur1} in Lemma \ref{lemma:strongly_convex} over $k$ and upper bounding $\rho_{k}$ with $\rho < 1$ we get:

\begin{subequations}
\begin{align}
    &\sum_{i=1}^{k} \norm{x_{i}-x_{i}^{*}} \leq \sum_{i=1}^{k} \rho \norm{x_{i-1}-x_{i}^{*}} + \sum_{i=1}^{k}\alpha \norm{e_{i}} + \sum_{i=1}^{k} \epsilon_{i}
    \\ 
    & \leq \rho \sum_{i=1}^{k} \norm{x_{i-1}-x_{i}^{*}} + \alpha E_{k} + P_{k}
    \\
    & \leq \rho \sum_{i=1}^{k} \norm{x_{i-1}-x_{i-1}^{*}} + \rho \sum_{i=1}^{k} \norm{x_{i-1}^{*}-x_{i}^{*}} + \alpha E_{k} + P_{k}
    \\
    & \leq \rho \sum_{i=1}^{k} \norm{x_{i}-x_{i}^{*}} + \rho \norm{x_{0}-x_{0}^{*}} + \rho \Sigma_{k} + \alpha E_{k} + P_{k} \label{last}.
\end{align}
\end{subequations}
Next, moving the first term on the right hand side of \eqref{last} to the left, we get~\eqref{eq:strongly_convex}.

\subsection{Proof of Theorem~\ref{thm:strongly_convex_maximum}}
\label{sec:proof_stronglyconvex_maximum}
Assume there exist finite constants $\sigma$, $\gamma_{e}$, $\gamma_{\epsilon}$ and $\rho$ such that $\sigma_k \leq \sigma$, $\norm{e_k} \leq \gamma_{e}$, $\epsilon_k \leq \gamma_{\epsilon}$  and $\rho := \sup_k \{\rho_k\} < 1$ for all $k \in \cT$. Then we can write \eqref{after_recur1-1} in Lemma \ref{lemma:strongly_convex} as follows:
\begin{equation}
    \norm{x_{k}-x_{k}^{*}} \leq \rho^{k} \norm{x_{0}-x_{0}^{*}}+ \frac{1-\rho^{k}}{1-\rho} (\rho \sigma + \alpha \gamma_{e} + \gamma_{\epsilon}) \, .
\end{equation}
Since $\rho <1$, as $k \to \infty$ we get~\eqref{eq:strongly_convex2}.

\subsection{Proof of Theorem~\ref{thm:regret_strongly_convex}}
\label{sec:proof_regret_strongly_convex}
Suppose that Assumptions~\ref{as:function_g}--\ref{as:minimum_attained} hold  and assume that there exists $D_k < + \infty$ such that $\|\partial f_k\| \leq D_k$ over the set $\cX_k$. Also define $D:= \sup_{k} \{D_{k}\}$. Then the dynamic regret can be written as:
\begin{subequations}
\begin{align}
    & \sum_{i=1}^{k} f_{i}(x_{i}) - f_{i}(x_{i}^{*}) \leq \sum_{i=1}^{k} 
 \langle \partial f_{i}(x_{i}), x_{i}-x_{i}^{*} \rangle
    \label{innerpro_thm3}
    \\
    & \leq \sum_{i=1}^{k} \norm{\partial f_{i}(x_{i})}\norm{x_{i}-x_{i}^{*}} 
    \\
    & \leq \frac{D}{1-\rho} (\rho \norm{x_{0}-x_{0}^{*}} + \rho \Sigma_{k} + P_{k} + \alpha E_{k})
     \label{cauchy_thm3}
\end{align}
\end{subequations}
where we have used Cauchy-Schwartz inequality and \eqref{eq:strongly_convex} in Theorem \ref{thm:strongly_convex} to go from \eqref{innerpro_thm3} to \eqref{cauchy_thm3}. Therefore, the result follows. 

\subsection{Proof of Theorem~\ref{thm:regret_convex}}
\label{sec:proof_regret_convex}
Since $g_{k}$ has a $L_{k}$-Lipschitz continuous gradient:
\begin{equation}
\label{Lg1}
    \begin{array}{l}
        g_{k}(x_{k}) \leq g_{k}(x_{k-1}) + \langle \nabla g_{k}(x_{k-1}),x_{k}-x_{k-1} \rangle 
        \\
        \hspace{3.6cm}+ \frac{L_{k}}{2} \norm{x_{k}-x_{k-1}}^{2} \, .
    \end{array}
\end{equation}
Using the convexity of $g_{k}$ we also get
\begin{equation}
\label{Cg1}
    \begin{array}{l}
        g_{k}(x_{k-1}) \leq g_{k}(x_{k}^{*}) + \langle \nabla g_{k}(x_{k-1}) , x_{k-1}-x_{k}^{*}\rangle.
    \end{array}
\end{equation}
Therefore, putting \eqref{Lg1} and \eqref{Cg1} together:
\begin{equation}
\label{Lg1&Cg1}
    \begin{array}{l}
         g_{k}(x_{k}) \leq g_{k}(x_{k}^{*}) + \langle \nabla g_{k}(x_{k-1}) , x_{k}-x_{k}^{*}\rangle 
         \\
         \hspace{4.5cm} + \frac{L_{k}}{2} \norm{x_{k}-x_{k-1}}^{2}.
    \end{array}
\end{equation}
On the other hand, rewrite equation \eqref{for_NSC_proof} as
\begin{equation}
\label{for_NSC_proof_rewrite}
\begin{array}{l}
    h_{k}(x_{k}) \leq h_{k}(x_{k}^{*}) - \frac{1}{\alpha} \langle x_{k-1}-x_{k}+r_{k},x_{k}^{*}-x_{k} \rangle 
    \\
    \hspace{2.2cm} + \langle \nabla g_{k}(x_{k-1})+e_{k},x_{k}^{*}-x_{k} \rangle + \frac{\epsilon_{k}^{2}}{2\alpha}
\end{array}
\end{equation}
with $\norm{r_{k}} \leq \epsilon_{k}$. Adding the  inequalities \eqref{Lg1&Cg1} and \eqref{for_NSC_proof_rewrite}:
\begin{equation}
    \begin{array}{l}
        g_{k}(x_{k}) + h_{k}(x_{k}) \leq g_{k}(x_{k}^{*}) + h_{k}(x_{k}^{*}) + \frac{L_{k}}{2} \norm{x_{k}-x_{k-1}}^{2} 
        \\ \hspace{1cm} + \frac{\epsilon_{k}^{2}}{2\alpha} - \frac{1}{\alpha} \langle  x_{k-1}-x_{k}+r_{k} , x_{k}^{*}-x_{k} \rangle + \langle e_{k}, x_{k}^{*}-x_{k}\rangle
        \\ \hspace{1cm} + \langle \nabla g_{k}(x_{k-1}) , x_{k}-x_{k}^{*}\rangle - \langle \nabla g_{k}(x_{k-1}) , x_{k}-x_{k}^{*}\rangle.
        \label{add_sub_thm4}
    \end{array}
\end{equation}
Adding and subtracting $x_{k}^{*}$ in the last term in the first line, and also in the second term in the second line of \eqref{add_sub_thm4} followed by using Cauchy-Schwartz inequality we get:
\begin{equation}
    \begin{array}{l}
        f_{k}(x_{k}) \leq f_{k}(x_{k}^{*}) + (\frac{L_{k}}{2}-\frac{1}{\alpha}) \norm{x_{k}-x_{k}^{*}}^{2} + \frac{L_k}{2} \norm{x_{k-1}-x_{k}^{*}}^{2} 
        \\
        \hspace{1cm} + \abs{\frac{1}{\alpha}-L_{k}}\norm{x_{k}-x_{k}^{*}} \norm{x_{k-1}-x_{k}^{*}}  + \frac{\epsilon_{k}^{2}}{2\alpha}
        \\
        \hspace{1cm} - \langle \frac{1}{\alpha} r_{k} - e_{k} , x_{k}^{*}-x_{k}\rangle \, .
    \end{array}
    \label{before_upperbdd}
\end{equation}

Set $\alpha \leq 1/\sup\{L_k\}$, say $\alpha = 1/\sup\{L_k\} - c^2$; then:
\begin{multline}
\frac{L_{k}}{2}-\frac{1}{\alpha} \leq \frac{\sup\{L_{k}\}}{2}-\frac{1}{\alpha} \leq   \frac{\alpha\sup\{L_{k}\} - 2}{2\alpha}  \leq \\ -\frac{1}{2 \alpha} - \frac{c^2 \sup\{L_{k}\}}{2 \alpha}.
\label{L_upperbdd}
\end{multline}

Also notice that $(\frac{1}{\alpha}-L_{k})\norm{x_{k}-x_{k}^{*}} \norm{x_{k-1}-x_{k}^{*}} \leq (\frac{1}{\alpha}-L_{k})\norm{x_{k}-x_{k}^{*}} (\norm{x_{k-1}-x_{k-1}^{*}} + \sigma_k)$. Let $R_k$ be  the diameter of $\cX_k$, and let $R$ be an upper bound on $\{R_k\}$. Therefore
\begin{equation}
    (\frac{1}{\alpha}-L_{k})\norm{x_{k}-x_{k}^{*}} (\norm{x_{k-1}-x_{k-1}^{*}} + \sigma_k) \leq \beta R (R + \sigma_k),
    \label{Diam_upperbdd}
\end{equation} 

where $\beta := \frac{1}{\alpha}- \inf\{L_{k}\}$ for brevity. Then, based on \eqref{L_upperbdd} and \eqref{Diam_upperbdd}, and neglecting constant negative terms, one can write \eqref{before_upperbdd} as:
\begin{equation}
    \begin{array}{l}
        f_{k}(x_{k}) \leq f_{k}(x_{k}^{*}) - \frac{1}{2\alpha} \norm{x_{k}-x_{k}^{*}}^{2} + \frac{1}{2\alpha} \norm{x_{k-1}-x_{k}^{*}}^{2}
        \\ \hspace{.1cm} + \frac{\epsilon_{k}^{2}}{2\alpha} - \langle \frac{1}{\alpha} r_{k}-e_{k}, x_{k}^{*}-x_{k}\rangle + \beta R (R + \sigma_k). \label{eq f}
    \end{array}
\end{equation}
Since $\norm{r_{k}} \leq \epsilon_{k}$, we can write \eqref{eq f} as follows:
\begin{equation}
\label{eq f1}
    \begin{array}{l}
        f_{k}(x_{k}) - f_{k}(x_{k}^{*}) \leq - \frac{1}{2\alpha} \norm{x_{k}-x_{k}^{*}}^{2} + \frac{1}{2\alpha} \norm{x_{k-1}-x_{k}^{*}}^{2}
        \\ + \frac{\epsilon_{k}^{2}}{2\alpha} + (\frac{\epsilon_{k}}{\alpha} + \norm{e_{k}}) \norm{x_{k}^{*}-x_{k}} + \beta R (R + \sigma_k).
    \end{array}
\end{equation}
Adding and subtracting $x_{k-1}^{*}$ in the second term on the right hand side of \eqref{eq f1}, and summing it from $i = 1$ to $k$, we can write the first two terms on the right hand side of \eqref{eq f1} as:
\begin{equation}
\label{eq f2}
    \begin{array}{l}
         \displaystyle \sum_{i=1}^{k} \textstyle \Big\{ \frac{-1}{2\alpha}\norm{x_{i}-x_{i}^{*}}^{2} + \frac{1}{2\alpha}\norm{x_{i-1}-x_{i-1}^{*}+x_{i-1}^{*}-x_{i}^{*}}^{2} \Big\}.
    \end{array}
\end{equation}
Then, we can upper bound \eqref{eq f2} as follows:
\begin{equation}
\label{eq f3}
    \begin{array}{l}
         \displaystyle \sum_{i=1}^{k} \textstyle \Big\{ \frac{-1}{2\alpha}\norm{x_{i}-x_{i}^{*}}^{2} + \frac{1}{2\alpha} \norm{x_{i-1}-x_{i-1}^{*}}^{2} 
        \\ + \frac{1}{2\alpha} \norm{x_{i-1}^{*}-x_{i}^{*}}^{2} +\frac{1}{\alpha} \norm{x_{i-1}-x_{i-1}^{*}}\norm{x_{i-1}^{*}-x_{i}^{*}} \Big\}.
    \end{array}
\end{equation}
The first two terms in \eqref{eq f3} form a telescoping cancellation; therefore, based on \eqref{eq f3}, and by considering other terms in \eqref{eq f1}, we can upper bound \eqref{eq f1} as follows: 
\begin{subequations}
\label{ave_reg}
\begin{align}
        & \sum_{i=1}^{k}  \big[ f_{i}(x_{i}) - f_{i}(x_{i}^{*})\big] \leq - \frac{1}{2\alpha} \norm{x_{k}-x_{k}^{*}}^{2} 
        \\ 
        & \hspace{1.0cm}+ \frac{1}{2\alpha} \norm{x_{0}-x_{0}^{*}}^{2}+ \frac{1}{2\alpha}\displaystyle \sum_{i=1}^{k} \sigma_{i}^{2} + \frac{1}{2\alpha} \displaystyle \sum_{i=1}^{k} \epsilon_{i}^{2}
        \\ 
        & \hspace{1.0cm}+ \displaystyle \sum_{i=1}^{k}  \sigma_{i} (\frac{1}{\alpha}\norm{x_{i-1}-x_{i-1}^{*}} +\beta R) + k\beta R^{2}
        \\ 
        & \hspace{1.0cm} + \frac{1}{\alpha} \displaystyle \sum_{i=1}^{k}  \big[(\epsilon_{i} + \alpha \norm{e_{i}}) \norm{x_{i}^{*}-x_{i}}\big].
\end{align}
\end{subequations}
Since $\cX_k \subseteq \cX$ is compact for all $k \in \cT$, we can upper bound $\norm{x_{k}-x_{k}^{*}}$ by $R$ and the result follows.

\bibliographystyle{IEEEtran}
\bibliography{ref.bib}

\end{document}